\documentclass[12pt, reqno]{amsart}
\usepackage{ifthen,amsfonts,amsmath,amssymb,epic,eepic}
\usepackage{epsfig}

\makeatletter
\@addtoreset{equation}{section}
\makeatother

\theoremstyle{definition}

\def\fnum{equation}
\newtheorem{Thm}[\fnum]{Theorem}
\newtheorem{Cor}[\fnum]{Corollary}

\newtheorem{Lem}[\fnum]{Lemma}

\newtheorem{Rem}[\fnum]{Remark}

\numberwithin{equation}{section}

\newcommand{\nn}{{\bf{n}}}
\newcommand{\Ric}{{\text{Ric}}}

\newcommand{\dist}{{\text {dist}}}

\newcommand{\mint}{-\!\!\!\!\!\!\int}

\def\RR{{\bf  R}}
\def\ZZ{{\bf  Z}}
\def\SS{{\bf  S}}

\newcommand{\e}{{\text {e}}}
\newcommand{\Area}{{\text {Area}}}

\newcommand{\cG}{{\mathcal{G}}}

\newcommand{\cL}{{\mathcal{L}}}

\newcommand{\eqr}[1]{(\ref{#1})}

\begin{document}

\title[Metrics without Morse index bounds]
{Metrics without Morse index bounds}

\author{Tobias H. Colding}%
\address{Courant Institute of Mathematical Sciences\\
251 Mercer Street, New York, NY 10012\\
MIT, 77 Massachusetts Avenue, Cambridge, MA 02139-4307\\
and Princeton University, Fine Hall, Washington Rd.,
Princeton, NJ 08544-1000}
\author{Nancy Hingston}%
\address{Department of Mathematics\\
The College of New Jersey\\
Ewing, NJ 08628.}

\thanks{The first author was partially supported by NSF Grant DMS 9803253.}

\email{colding@cims.nyu.edu and hingston@TCNJ.EDU}

\maketitle

\section{Introduction}
Let $M^2$ be a closed orientable surface with curvature $K$
and $\gamma\subset M$
a closed geodesic.  The {\it Morse index} of
$\gamma$ is the index of the
critical point $\gamma$ for the length functional, i.e., the number of
negative
eigenvalues (counted with multiplicity) of the second
derivative of length (throughout curves will always be in $H^1$).  
Since the second derivative of length at $\gamma$
in the direction of a normal variation $u\,\nn$ is
$-\int_{\gamma}u\,L_{\gamma}\,u$ where $L_{\gamma} \,u= u'' + K\,u$,
the Morse index is the number of
negative eigenvalues of $L_{\gamma}$.
(By convention,  an
eigenfunction $\phi$ with eigenvalue
$\lambda$ of $L_{\gamma}$ is a solution of
$L_{\gamma}\,\phi+\lambda\, \phi=0$.)  Note that if $\lambda=0$,
then $\phi$ (or $\phi\,\nn$) is a (normal) Jacobi field.
$\gamma$ is {\it stable}
if the index is zero.
The {\it index} of a noncompact geodesic is the dimension
of a maximal vector space of compactly supported variations for which the
second derivative of length is negative definite.  We also say that
such a geodesic
is {\it stable} if the index is $0$.  

Our main result is:

\begin{Thm}  \label{t:example}
On any $M^2$,
there exists a metric with a geodesic lamination
with infinitely many unstable leaves.  Moreover, there is such a metric
with simple closed geodesics of arbitrary high Morse index.
\end{Thm}

The first part of Theorem \ref{t:example} is relatively easy to 
achieve and in the proof we do that first.  

A codimension one {\it lamination}
on a surface $M^2$ is a collection $\cL$ of
smooth disjoint curves (called leaves)
such that
$\cup_{\ell \in \cL} \ell$ is closed.
Moreover, for each $x\in M$ there exists an
open neighborhood $U$ of $x$ and a $C^0$ coordinate chart, $(U,\Phi)$, with
$\Phi (U)\subset \RR^2$
so that in these coordinates the leaves in $\cL$
pass through $\Phi (U)$ in slices of the
  form $(\RR\times \{ t\})\cap \Phi(U)$.  
A {\it foliation} is a lamination for which 
the union of the leaves is all of $M$
and a {\it geodesic lamination} is a lamination whose leaves are geodesics.

Theorem \ref{t:example} will be proven by first constructing a 
metric on the disk with convex boundary having no Morse index bounds 
and then completing the metric to a metric on the given  
$M^2$.  By taking the product of this metric on the disk 
with a circle we get, on a solid torus, 
a metric with convex boundary and without Morse index bounds for 
embedded minimal tori, and with a minimal lamination with infinitely 
many unstable leaves.  By completing this metric we get: 

\begin{Thm}  \label{t:example2}
On any $M^3$,
there exists a metric with a minimal lamination
with infinitely many unstable leaves.  Moreover, there is such a metric
with embedded minimal tori of arbitrary high Morse index.
\end{Thm}

By construction the embedded minimal tori in Theorem \ref{t:example2} and the 
leaves of the lamination can be taken to be totally geodesic.

Our interest in whether the Morse index is bounded comes from the 
assertion of J. Pitts and J.H. Rubinstein (see \cite{PiRu}, \cite{CM2}) 
that if one can show that the Morse index of all embedded
minimal tori is bounded for a sufficiently large class of metrics on
$\SS^3$, then the
spherical space form problem can be solved affirmatively.

We will equip the space of metrics on a given manifold with the
$C^{\infty}$-topology.  A subset of the set of metrics on the manifold is said
to be {\it residual} if it is a countable intersection of open dense
subsets.
A metric on a surface is {\it bumpy} if each  closed geodesic
is a nondegenerate critical point, i.e.,
$L_{\gamma} u = 0$ implies $u\equiv 0$.  Bumpy metrics are generic,
\cite{Ab}, \cite{An}; that is the set of bumpy metrics contain a
residual set.

To check that any given metric is bumpy is virtually impossible; however 
it seems that the metric in Theorem \ref{t:example} can be chosen 
to be bumpy; see also Remark \ref{r:genright} below.  Thus it seems 
unlikely that a bumpy metric is enough to ensure a bound for the Morse 
index of simple closed geodesics on $M^2$.  What is needed is a 
nondegeneracy condition for noncompact simple geodesics, rather 
than one for closed geodesics; cf. \cite{CH1}, \cite{CH2}.

Bounding the Morse index can be thought of
as a purely
analytical problem about lower bounds for
eigenvalues of a Schr\"odinger
operator.  A typical way of getting such lower bounds is in terms of 
integral bounds for the potential.  For instance, if $L\, u = u'' + K u$
is a Schr\"odinger operator on a
circle ${\mathcal{C}}$ with
$\text{Length}({\mathcal{C}})\int_{{\mathcal{C}}} \max \{ K , 0 \}
\leq C$, then
the number of negative eigenvalues counted with multiplicity
is bounded by $N = N(C)$.
However, as the next theorem illustrates,
bounding the Morse index in this setting 
is analytically rather subtle; see also \cite{CH1}, \cite{CH2}.

\begin{Thm}  \label{t:e0.1}
On any $M^2$,
there exists an open (nonempty) set of metrics
so that for each metric there is
a sequence of simple closed stable geodesics $\gamma_i$ with
$\text{Length}(\gamma_i)\to \infty$ and
\begin{equation}   \label{e:e02}
\inf_i  \mint_{\gamma_i}
\max \{K,0\} \equiv
\inf_i    \int_{\gamma_i}
\max \{K,0\} \, / \,  \text{Length}(\gamma_i)
 >0\, .
\end{equation}
\end{Thm}

Similarly,
for $3$-manifolds and minimal surfaces (with second fundamental
form $A$):

\begin{Thm}  \label{t:e0.1a}
On any $M^3$,
there exists an open (nonempty) set of metrics
so that for each metric there
is a sequence of embedded stable minimal tori $\Sigma_i$ with
$\Area (\Sigma_i)\to \infty$ and
\begin{equation}   \label{e:e02b}
\inf_i\mint_{\Sigma_i}\max \{|A|^2+\Ric_M (\nn,\nn),0\}>0\, .
\end{equation}
\end{Thm}

Theorems \ref{t:e0.1}, 
\ref{t:e0.1a} are 
weaker than Theorems \ref{t:example}, \ref{t:example2} 
in the sense that they do 
not produce examples of metrics with no index bound.  But the open set of 
metrics given in Theorems \ref{t:e0.1}, \ref{t:e0.1a} means that it is 
impossible to prove that 
the Morse index is bounded for a generic metric using only the standard 
analytic argument mentioned just above Theorem \ref{t:e0.1}.  
It follows that a bound for the Morse bound is not just a simple analytical 
fact but relies on the
(global) dynamics of the situation.  This is where a generic condition is 
needed.

Since the metrics of 
Theorems \ref{t:e0.1} and 
\ref{t:e0.1a} are much simpler than those of Theorems \ref{t:example} 
and \ref{t:example2}, then they are described first in Section 
\ref{s:schr}.

\vskip2mm
Throughout this paper (except in Theorems \ref{t:example2}, 
\ref{t:e0.1a})
$M^2$ is a closed orientable surface with a Riemannian metric, $\cL$
is a geodesic lamination,
and when $x\in M$, $r_0>0$, and
$D\subset M$, then we let $B_{r_0}(x)$ denote the ball of
radius $r_0$ centered at $x$ and $T_{r_0}(D)$
the $r_0$-tubular neighborhood of $D$.  Moreover, if $x,$ $y\in M$, then
$\gamma_{x,y}:[0,\dist_M(x,y)]\to M$ will denote a minimal geodesic
from $x$ to $y$.  Whenever we look at a single geodesic it will always
be assumed to be parameterized by arclength.

\vskip2mm
We are grateful to Nicholas Hingston Tenev for making the figures.

\section{Index of Schr\"odinger operators}   \label{s:schr}

In this section, we will discuss some of the difficulties with approaching
the problem of bounding the Morse index from a purely analytical point of
view.  In particular, we will show Theorems \ref{t:e0.1} and \ref{t:e0.1a}.

\begin{proof}
(of Theorem \ref{t:e0.1}).
Let $S \subset M$ be a connected planar domain with
three interior boundary components and $\nu:[0,1]\to S$
a simple curve connecting two different components of $\partial S$.

It is easy to see from the proof of
proposition $1.1$ of \cite{CM1} that it is enough to show
that the set of metrics
$\cG$ on $S$ so that $S$ has strictly convex boundary and $\min_{\nu}K>0$ is
nonempty (observe that $\cG$ can clearly be extended to an open set of
  metrics on all of $M$).  This will ensure that the curve $\gamma_m$ in 
Figure 1 will cross a region of positive 
curvature $m$ times, if $\nu$ connects the two lower interior 
boundary components.

To see that $\cG$ is nonempty, let $\sigma_1$ and $\sigma_2$ be the two
boundary components of $S$ that $\nu$ connects.  It is clearly enough to
show that we can find a metric on a neighborhood (in $S$) of
$\sigma_1\cup \sigma_2\cup \nu$ such that $\min_{\nu}K>0$ and the
outward (to $S$)
normal geodesic curvatures of $\sigma_1$ and $\sigma_2$ are positive.
This can easily be achieved by letting $B_{r_0}$ be a ball of radius
$r_0<\pi/2$ on the unit
sphere $\SS^2$ and $\nu\subset B_{r_0}$ a
geodesic through the center of the ball and thinking of $\sigma_1$ and
$\sigma_2$ as two disjoint copies of $\partial B_{r_0}$ each of
which intersects $\nu$ in only one end.
\end{proof}

\begin{figure}
    \setlength{\captionindent}{4pt}
    \begin{minipage}[t]{0.5\textwidth}
    \centering\includegraphics{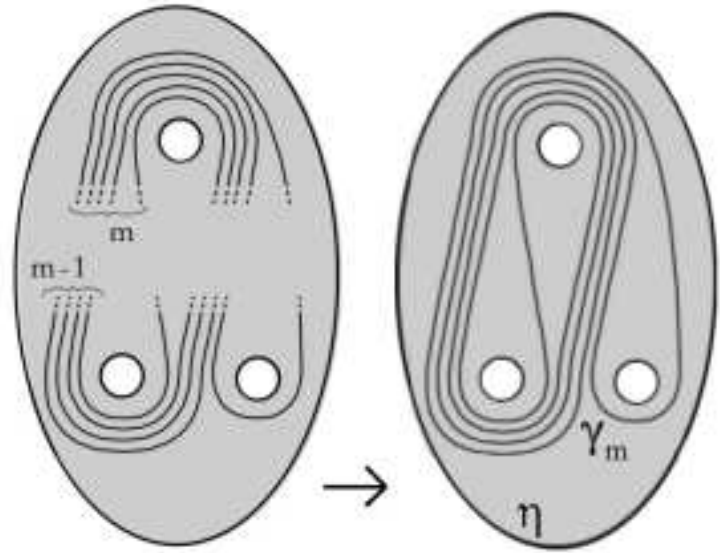}
    \caption{The curves $\gamma_m$ on a disk with 3 holes.}
    \end{minipage}%
\end{figure}

Observe that if $\gamma\subset M^2$
is a closed stable geodesic, then by the stability inequality
applied to the function $\phi\equiv 1$
\begin{equation} \label{e:e03}
\int_{\gamma}K \leq 0\, .
\end{equation}
Likewise if $\Sigma\subset M^3$ is a closed stable minimal surface, then
\begin{equation}
\int_{\Sigma}[|A|^2+\Ric_M (\nn,\nn)] \leq 0\, .
\end{equation}

On the other hand, in one dimension (and similarly in two
dimensions; cf. Theorem \ref{t:e0.1a}),
easy examples show that if the potentials $K_i$ of a sequence of
Schr\"odinger
operators on circles ${\mathcal{C}}_i$ satisfy
$
\text{Length}({\mathcal{C}}_i)\int_{{\mathcal{C}}_i}\max\{K_i,0\}\to \infty
$
and
$\int_{{\mathcal{C}}_i} K_i \leq 0$, then typically
there is no uniform bound for the indices.

\begin{proof}
(of Theorem \ref{t:e0.1a}).
Let again $S$ be a connected planar domain with
three interior boundary components and $\nu:[0,1]\to S$
a simple curve connecting two different components
$\sigma_1$ and $\sigma_2$ of
$\partial S$.  Set $\Omega=S\times \SS^1$.  By the proof of
Theorem $0.1$ of \cite{CM1}, it
is clearly enough to show that there exists a metric on $\Omega$ with
strictly convex boundary and such that $\min_{\nu\times\SS^1}K_{\Omega}>0$.
In fact it is enough to show that such a metric exists in a neighborhood
of $\sigma_1\times\SS^1\cup\sigma_2\times \SS^1\cup\nu\times\SS^1$
(in $\Omega$).

To see this observe first that clearly there exists such a
metric in a neighborhood of $\nu\times\SS^1$ (in $\Omega$) - we need to
extend this to a metric in a neighborhood of all of
$\sigma_1\times\SS^1\cup\sigma_2\times \SS^1\cup\nu\times\SS^1$
while keeping the surfaces $\sigma_1\times\SS^1$ and $\sigma_2\times\SS^1$
strictly convex.  This can easily be achieved by first choosing any
extension to $\sigma_1\times \SS^1\cup \sigma_2\times \SS^1$ and then
extending it to a small normal neighborhood of
$\sigma_1\times \SS^1\cup \sigma_2\times \SS^1$ making
$\sigma_1\times \SS^1\cup \sigma_2\times \SS^1$ strictly convex.
\end{proof}

\section{Metrics without Morse index bounds}   \label{s:s5}

In this section we will prove Theorems \ref{t:example} and \ref{t:example2}.
As mentioned in the introduction it suffices to construct a 
metric on the disk 
with convex boundary that has a geodesic lamination with infinitely many 
unstable leaves, and show that there are  simple
closed geodesics of arbitrarily high Morse index. The
construction
relies on the:

{\sc Basic Barrier Principle:}  Let $R$ be a domain in
$M^2$ with piecewise geodesic boundary that is locally
convex in $M$; that is, $R$ should locally be the intersection of
closed geodesic half-spaces in $M$. Given a simple closed
curve $\gamma$ on $R$, there is a simple closed
geodesic $\tau$ which is freely homotopic to $\gamma$ on $R$,
and which has length less than or equal to that of $\gamma$.

The boundary of $R$ should be thought of as a
barrier.  The principle is a consequence of a
``Lusternik-Schnirelmann curve-shortening process'';
that is, of a
curve shortening process which does not introduce
intersections. There are many versions of such a
process; see for instance
Grayson's paper \cite{Gr}, and the related work of
Angenent \cite{Ang}. With this process a closed curve
evolves in the direction of its curvature vector
$\kappa \,\nn$.  Grayson shows that, if the curve does
not
shrink to a point in finite time, then the flow is
defined for all time, and that for every $n\geq 0$,
\begin{equation}  \label{e:kappainfty}
\lim_{t\rightarrow \infty} \sup \kappa^{(n)}=0\, .
\end{equation}
Angenent has shown that the number of intersection
points of two evolving curves can only decrease.  Using
this,
we get a slightly stronger principle:

{\sc Second Barrier Principle:} Let $\tilde{M}$ be a
covering space of $M$, and let $\tilde{R}$ be a locally
convex
domain in $\tilde M$ as before. (We are not assuming
that
$\tilde{R}$ is the lift of a set $R$.)
Let $\tilde \gamma$ be a simple closed curve on $\tilde
{R}$ whose
image on $M$ is still simple.  Then there is
a simple closed geodesic $\tilde \tau$ in the same
free homotopy class as $\tilde\gamma$ on $\tilde{R}$
whose image $\tau$
on $M$ is still simple, and which has length
less than or equal to that of $\tilde{ \gamma}$.
Moreover, if $\sigma$ is a
geodesic on $\tilde{R}$ (closed, or with boundary on
$\partial \tilde{R}$), then $\tilde{\tau}$ intersects
$\sigma$ in no
more points than did  $\tilde{\gamma}$.  If the metric
is bumpy, then we can
assume that $\tilde{\tau}$ is a local
minimum of length.  If there is no geodesic in $\partial
\tilde{R}$ homotopic to $\tilde{\gamma}$, then $\tilde \tau$ 
will lie in the interior of $\tilde R$.

\begin{proof}
This follows since there are no choices involved
in the curve-shortening process.
\end{proof}

   We start with a surface of
revolution about the
$z$-axis which is the connected sum of three round
spheres
using two (sufficiently narrow) necks of negative
curvature.  It looks like a
   snowman.  It will be useful below to refer
to the region between the
necks as the ``middle-sphere'' part of
$M$.
Let $(r,\theta, z)$ be cylindrical
coordinates, and let $\alpha$ be an angle defined on
the tangent bundle to the surface of revolution as
the angle measured down clockwise from the meridian
$\theta$= constant.
(So $\alpha =0$ for a vector on
the surface which points up; $\alpha =\pi /2$ for a
vector which points in the direction of increasing
$\theta$.)
Curves with $dz/d\theta>0$ will be called
{\em positive}.  Later we will alter the metric on our surface.  
The coordinate functions $\theta$, $z$, and $\alpha$ on the surface 
will be fixed however, induced from the embedding as a surface of 
revolution in $\RR^3$.

The relative radii of the three balls are not
important, but the
two necks must have equal
radii; this way by Clairaut's theorem there will be
two 1-parameter
families of geodesics, which we will call
{\em limiting geodesics}, crossing the (middle-sphere)
equator with angles $\pm \alpha_{0}$, which
spiral toward the upper neck as $t\rightarrow \infty$
and toward the lower neck as $t\rightarrow -\infty$.
(Geodesics crossing the equator with $| \alpha
|\leq\alpha_{0}$
will enter the lower and upper spheres;
those with $| \alpha |  > \alpha_{0}$ will stay in the
center sphere with periodic $z$-values.)
It will be convenient to have each neck symmetric about
a horizontal plane; this way a geodesic which crosses a
neck will emerge as an inverted mirror image of the
curve that went in.

Next we add a nose to the top sphere and two to
the bottom sphere. For each nose we alter the metric
in a small disk away from the $z$-axis in such a way
that the disk now supports a strictly stable simple
closed
geodesic which
we will use as a barrier.  (The shape of the nose is
unimportant so long as we have this barrier geodesic.
These barrier
geodesics will serve much the same purpose as the
interior boundary of
$S$ did in the proof of Theorem \ref{t:e0.1}.)
It will also be useful to have two more simple strictly
stable
geodesic
barriers circling the $z$-axis
at the top and bottom of the snowman.  To get a disk, 
we cut away the interior of the lower barrier circling the $z$-axis.

We call this Riemannian surface $M_{0}^{2}$.  (See
Figure 2.)

\begin{figure}
    \setlength{\captionindent}{4pt}
    \begin{minipage}[t]{0.5\textwidth}
    \centering\includegraphics{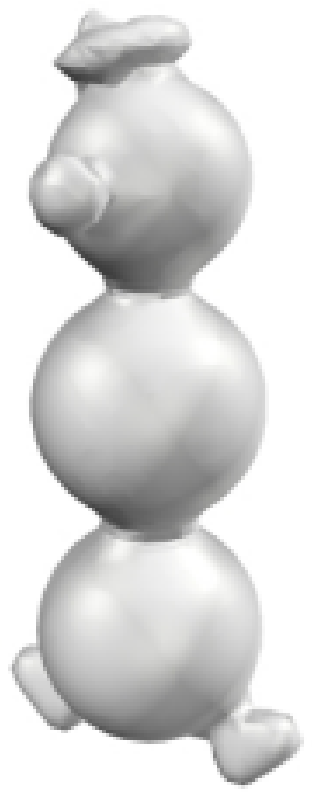}
    \caption{$M_0^2$.}
    \end{minipage}%
\end{figure}

We need one more alteration of the metric:  We will add
an
infinite number of ``bumps''  along the equator.  
Given a radius $\rho$ and a point $P$ on the equator of
the
middle-sphere, we will make a ``bump'' of radius $\rho$
at $P$ with the following four properties:\\
    1. The metric remains unchanged except within
distance $\rho /2$
    of $P$.\\
    2. The metric remains locally rotationally symmetric
about
    $P$.\\
    3. In the new metric, the geodesics through $P$ are
    longer than they were without the bump.\\
    4. The metric remains smooth after adding an infinite
    number of disjoint bumps.\\
    This can be done for example by multiplying the metric by a 
conformal factor $1+\e^{-1/\rho}f(r/\rho)$ in coordinates centered 
at $P$ ($r$ is the radial component) for some smooth bump
function $f$ with support in $[0,1/2]$.  In order to prove the first 
statement in Theorem \ref{t:example}, it will be convenient 
to have $f$ monotonic.  We place the bumps 
around the equator in such a way that each segment of length 
$\mu$ (to be determined below) contains an accumulation point of 
the bump centers.  They should be
far enough apart that between any two there is a
geodesic which crosses the equator at angle $\alpha_{0}$.

This surface with bumps is $M^2$.

 {\sc The Geodesics.}  We digress briefly to say
roughly what the geodesics of high index on $M^2$ look like.
First we describe some geodesics on the surface $M_{0}^{2}
$. Let $\eta$ be a simple closed curve with the three noses in 
its interior and the north- and south-pole barrier 
geodesics in its exterior.  We can assume by the  
barrier principle that $\eta$ is geodesic.  Let $R$ be the 
interior of $\eta$ with the three barrier holes removed. 
Figure 1 shows a simple closed curve $\gamma_{m}$ on 
$R$, which we can again assume is geodesic.  Here $m$ is the 
number of ``strands'', so that $\gamma_{m}$ crosses 
the  equator $2m$ times.  Now put a rubber band 
around $\eta$, and another around $\gamma_{m}$. Pick an 
integer $n> 0$. Take the snowman's head, and rotate 
it  (just the head) counterclockwise through an angle 
of $2n\pi$.  (Don't let the rubber bands get unhooked.)  
This results in new geodesics $\eta_{n}$ 
and $\gamma_{m,n}$ which wind (and unwind) $n$ times.  
Now these geodesics are still local minima (index $0$).  
But if $n$ is sufficiently large (depending upon $m$), we 
can use a minimax argument to produce a geodesic close to $\gamma_{m,n}$
    on the surface $M^{2}$ which goes over a bump
    each time it crosses the equator, and thus has index
$2m$.
    The bumps on the equator act as a comb which enables
us
    to separate the strands and move them around the
    equator.  Of course the bumps are arbitrarily small,
    and we need them high enough to hold the strands.
    However as $n\rightarrow \infty$ the strands of
$\gamma_{m,n}$
    approach the critical ``limiting''
    angle $\alpha _{0}$ and it takes less and less to
hold them.

\vskip2mm
We will now make all of this precise in the following
theorem:

    \begin{Thm} \label{t:counterex}
For every $k\geq 0$, the surface $M^2$ has a simple
closed geodesic with Morse index in
$[k-2,k]$.
\end{Thm}

The geodesics we will produce will have local homology
in dimension $k$, and thus (see \cite{K} corollary 2.5.6 
and proposition 1.12.3 and also \cite{Ch} theorem 5.4 on p. 50) index $\leq
k\leq$ index $+$ nullity $\leq$ index $+2$,
since on a surface any closed geodesic has
nullity $\leq 2$.

The proof of Theorem \ref{t:counterex} will rely on a
number of
observations and lemmas.
We begin by establishing some basic facts about the
curve $\eta_{n}$.  The curve $\eta_{n}$ is a curve on a
cylinder ($M_{0}^{2}$
    with north pole removed) with $3$
    holes removed.  It wraps around the cylinder $n$
times, hooks around the top nose, unwraps $n$ times, and
hooks around the bottom two noses.
    Using curve-shortening, we can assume that $\eta_n$ is geodesic, 
and that no other curve in the free homotopy class, and intersecting 
the equator exactly twice, has length less than that of $\eta_n$.  
Let the barrier circling the
north-pole lie at $z=z_{N}$, and the upper neck at
$z=z_{0}$.  Let $z_{1}> z_{0}$.  After possibly
making the necks more narrow, we can assume that there
is an $\epsilon > 0$ so that if $\sigma_{z}(\theta)$ is
the curve on $M_{0}^{2}$ at height $z$ parameterized by
$\theta$, then $|\sigma_{z}'(\theta)|\geq |\sigma_{z_
{0}}'(\theta)|+\epsilon$
for every $z$ with $z_{1}\leq z \leq z_{N}$.
Since $\eta_n$ has minimal length, for fixed $z_{1}$
the length of the part of $\eta_{n} $ which lies above
$z=z_{1}$ is bounded as $n\rightarrow\infty$.  
All the winding around takes place between $z=z_{2}
$ and $z=z_{1}$; thus the angles at which $\eta_{n}$ crosses the equator
both go to $\alpha_{0}$ as $n\rightarrow \infty$.

Let $R_{n}$ be the domain
    bounded by $\eta_{n}$ and the $3$ hole boundaries.
Let $E_{n}$ be the intersection of $R_{n}$ with the
equator.  We claim:
\begin{equation}
    \mu =\frac{1}{2}\liminf_{n\rightarrow\infty}\text{width}
(E_{n})>  0\, .
\end{equation}
To see this, note that the two sides of $\eta_{n}$
are bounded apart at the equator if and only if they are
    bounded apart as they cross a level
    $z=z_{1}$ just above the upper neck, and if they are
bounded apart as they cross a level $z=z_{2}$ just
below the lower neck.  Since they travel a bounded
distance outside the region $z_{2}\leq z\leq z_{1}$,
and since $\eta_{n}$ is simple,
this will be the case.

Since the set of bumpy metrics is dense (by \cite{Ab})
we can let ${g_{i}}$ be a sequence of bumpy metrics with limit
the surface $M^2$ as $i\rightarrow \infty$.  (We
introduce the metrics $g_{i}$ in order to simplify a
Morse theoretic argument which will come later.  We
will use Morse theory in the metrics $g_{i}$ to find
simple closed geodesics of high index in the metrics $g_
{i}$, where the length functional is a nondegenerate
Morse function, and then take the limit as
$i\rightarrow\infty$.) Using the
fact that the $5$ barrier geodesics are strictly
stable, we can (cf. Lemma \ref{l:exF}, 
Lemma \ref{l:vary} below)
assume that each metric $g_{i}$
has $5$ simple closed geodesics with limit as
$i\rightarrow\infty$ the $5$ barrier geodesics 
on $M^2$.  For $i$ sufficiently large, applying curve-shortening in $g_i$ 
metric to the curve $\eta_n$ will produce a simple closed geodesic 
$\eta_{i,n}$ with $\eta_{i,n}\to \eta_n$ as $i\to \infty$.  Here we are again 
using the fact that $\eta_n$ has minimal length; it might be necessary to 
replace $\eta_n$ by another simple closed geodesic of the same length, 
and with all the properties established above.  Note $\eta_{n}$ will still
be simple since on an orientable surface a geodesic
which is a limit of simple curves will itself be simple.

Next we show that the bumps do bend geodesics; a
geodesic crossing the equator to the right of a bump
center will
curve to the left at both ends.

\begin{Lem} \label{l:l5.1}
Consider the surface $M_{0}^{2}$ with one bump of
radius $\rho$ added at a point $P$ on the equator.  If
$\rho$ is sufficiently small, there is a geodesic
$\lambda$ which meets a ball of radius $\rho /2$ about $P$, and
crosses the equator to the right of $P$ at an angle
$\alpha \approx \alpha _{0}$.  The ``top half'' of
$\lambda$ will leave the ball of radius $\rho$ about $P$
at an angle more vertical than the critical angle, and
will cross each positive limiting geodesic an infinite
number of times in the northern part of the middle-sphere
before crossing the upper neck.  The ``bottom half''
of $\lambda$ will leave the ball of radius $\rho$
at an angle less vertical than the limiting angle, and
will cross each positive limiting geodesic at least
twice in the southern part of the middle-sphere before
crossing the equator again.  If $\tau$ is the limiting
geodesic which is tangent to the ball of radius $\rho$
about $P$ on its right hand side, then $\lambda$ stays
to the left of $\tau$ in the universal cover of the
cylinder.
\end{Lem}

\begin{proof}
Let $\zeta$ be the geodesic through $P$ at angle
$\alpha_{0}$.  This geodesic is longer than the
geodesic $\tau$ to its right which does not meet the
bump.  Since the two geodesics are asymptotic at both
ends, over a large enough distance $\zeta$ is not
length minimizing; over a large distance it becomes
more efficient to go around the bump.  Take two far
away points on $\zeta$, and find a minimizing geodesic
$\lambda$ ($\neq \zeta$) in the homotopy class determined by the 
segment $\zeta$.  We can
assume that $\lambda$ lies to the right of $\zeta$ and
to the left of $\tau$. Extending $\lambda$ at both ends
will produce a geodesic as described; after leaving the
bump, the geodesic $\lambda$ is described by Clairaut's
theorem.  (The hypothesis that $\rho$ is not too large
is used to ensure that the bottom half of $\lambda$ does
not cross the equator again too soon:  In the limit as
$\rho \to 0$ the bottom half of $\lambda$
becomes a positive limiting geodesic ``followed by'' a
negative limiting geodesic which will cross each
positive limiting geodesic an infinite number of times
before crossing the equator again.)
\end{proof}

\begin{Cor} \label{c:}
The surface $M^2$ has a geodesic lamination with
infinitely many unstable leaves.
\end{Cor}

\begin{proof}
The limiting geodesics in the metric $M_{0}^2$ which
pass through bump centers are still geodesics on $M^2$,
since the bumps are rotationally symmetric.  If $f$ is monotonic, 
it is easy to see that they have positive index. (See Figure 3.) 
The closure of their union is the desired lamination. 
\end{proof}

\begin{figure}
    \setlength{\captionindent}{4pt}
    \begin{minipage}[t]{0.5\textwidth}
    \centering\includegraphics{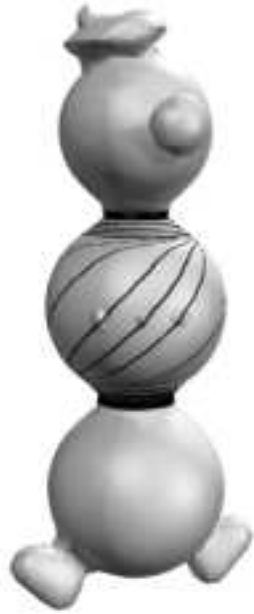}
    \caption{A geodesic lamination on $M$ 
with infinitely many unstable leaves.}
    \end{minipage}%
\end{figure}

This proves the first part of Theorem \ref{t:example}.

   We now return to the curves $\eta_{n}$.  Let $n$ and
$\rho$ be given, and fix a point $P$ on $E_{n}$. We
consider again the surface $M_{0}^{2}$ with a single
bump of radius $\rho$ added at $P$.   We say {\em the
bump of radius $\rho$ at $P$ is high enough to hold
$\eta_{n}$} if $1$) there is a geodesic segment
$\lambda$ on
$R_{n}$
with both endpoints on the left side of the
middle-sphere section of $\eta_{n}$ which crosses $E_{n}
$
exactly once, to the right of $P$; and if $2$) (the
mirror image of $1$)) there is a
geodesic segment $\psi$ on $R_{n}$ with both endpoints
on the
right side of the middle-sphere section of $\eta_{n}$
which crosses $E_{n}$ exactly once, to the left of
$P$.  (See Figure 4.) The idea will be to use $\lambda$
and $\eta_{n}$ together as barriers to get a strand of
$\gamma_{m,n}$ to hook around the right side of the
bump at $P$.

\begin{figure}
    \setlength{\captionindent}{4pt}
    \begin{minipage}[t]{0.5\textwidth}
    \centering\includegraphics{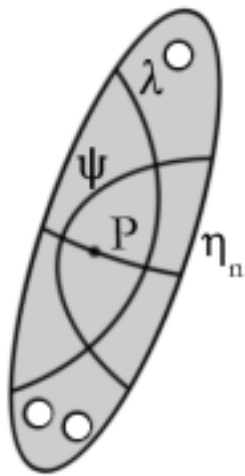}
    \caption{A bump at the point $P$ ``high enough to hold $\eta_n$'' deflecting the geodesics $\lambda$ and $\psi$.}
    \end{minipage}%
\end{figure}

\begin{Lem} \label{l:l5.2}
Fix $m\geq 1$.  If $n$ is sufficiently large, there are
$2m$ bumps on $R_{n}$ in $M^2$ which are high enough to
hold $\eta_{n}$.
\end{Lem}

\begin{proof}
By construction there is a $\epsilon>0$ so that each
interval on the equator of width $\mu$ contains at
least $2m$ bumps of radius at least $\epsilon$.  By
the definition of $\mu$, for $n$ sufficiently large $E_{n}$ 
will contain an interval of length $\mu$.  Now we
use Lemma \ref{l:l5.1} and the fact that the two
``sides" of
$\eta_{n}$ are geodesics at angles which approach the
limiting angle as $n\rightarrow\infty$. Of course Lemma
\ref{l:l5.1}
works as well to get a geodesic curving around the
left side of $P$.
\end{proof}

Next we bring in the curves $\gamma_{m,n}$. To reduce
the number of subscripts, fix $m$ and $n$ with $n$
``sufficiently large."  Here is the idea of the rest of the proof.  
For simplicity assume for the moment that the metric is bumpy.  
Let $P_{1},\cdots,P_{2m}$ be
points on $E_{n}\subset M^2$ which are bump centers for
bumps ``high enough to hold $\eta_{n}$," ordered from
left to right.  We will show that there is a simple
closed geodesic, freely homotopic to $\gamma_{m,n}$ in
$R_{n}$, where we can choose the $j$'th strand (ordered
left to right) to cross
the equator either to the left or the right of $P_{j}$.  
With a bumpy metric this geodesic will be a local minimum of length.  
This gives, with all possible such choices, $2^{2m}$ distinct simple 
closed geodesics of index $0$, each freely homotopic to $\gamma_{m,n}$.  
We next (Lemma \ref{l:l5.3}) fill in a $2m$-cube of simple closed curves 
freely homotopic to $\gamma_{m,n}$, with the $2^{2m}$ geodesics as vertices.  
We will use relative homology to show (Lemma \ref{l:l5.4}) 
that each face of the 
cube will ``lie hanging'' on a closed geodesic whose Morse index is equal 
to the dimension of the face.  (Note to get index $>1$ it is not 
sufficient to have a cube whose vertices are distinct local minima; 
to get nontrivial topology in higher dimensions we use \eqr{e:e5.12} below.)

\begin{figure}
    \setlength{\captionindent}{4pt}
    \begin{minipage}[t]{0.5\textwidth}
    \centering\includegraphics{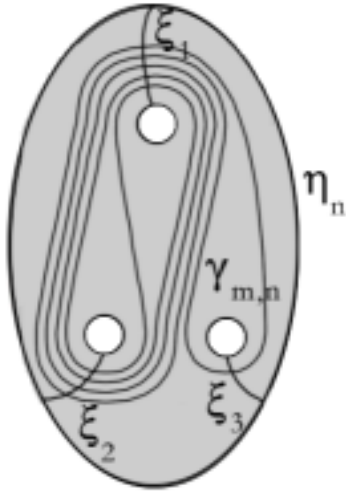}
    \caption{Geodesics in the domain $R_n$.}
    \end{minipage}%
\end{figure}

   Let
$\xi_{1}$, $\xi_{2}$, $\xi_{3}$ be simple curves in $R_
{n}$ starting on
the $3$
hole boundaries and ending on $\eta_{n}$.   They should
not intersect each
other or $E_{n}$. Figure 5 shows a domain homeomorphic
to $R_{n}$ (actually the homeomorphism type does not
depend on $n$.) Let $\gamma_{m,n}$ be the curve shown
in Figure 4.  We
assume that $\gamma_{m,n}$ has intersections with the
$\xi_{i}$ as in the
figure.  Let $S$ be the $2m$-fold
covering space of $R_{n}$ in which $\gamma_{m,n}$ lifts
to a closed curve $\tilde{\gamma}_{m,n}$ which crosses
a
different lift of
$E_{n}$ each time it crosses the equator. (The cover
can be constructed as
follows:  Get $2m$ copies of $R_{n}$ with the $\xi _{i}
$ marked.  Pick
one, and start tracing out the curve $\gamma _{m,n}$.
Each time you cross
one of the $\xi _{i}$, cut along $\xi _{i}$ and paste
in a new copy of
$R_{n}$.)  Let $\tilde{E_{n}}$ be the set in $S$ lying
above $E_{n}$. (It will have $2m$ components.)  At this
point it should be clear, using the barrier principle,
that we can assume that $\tilde{\gamma}_{m,n}$ is
geodesic and still has all the same intersections: It
crosses each component of $\tilde{E_{n}}$ once, and its
image on $M^2$ crosses each neck $2m$ times.  In
between
the necks its $z$-coordinate is monotone.
Label the intersection points of $\tilde{\gamma}_{m,n}$ with
$\tilde{E_{n}}$: $Q_{1},\cdots,Q_{2m}$ in such a way
that
their images in $E_{n}$ go from left to right.  Let ${T_
{j}}$ be the $2m$ lifts of the portion of
$R_{n}$ which lies between the two necks on $M^2$,
labeled so that $T_{j}$ contains $Q_{j}$. 

Given subsets $I$ and $J$ of $\{1, \ldots ,2m\}$, let
\begin{equation}
\ell_{j}=\left\{
\begin{array}{ll}
j & \mbox{if $j\in I$}\, ,\\ j-1 & \mbox{otherwise}\, ,
\end{array}
\right.
\end{equation}
\begin{equation}
r_{j}= \left\{
\begin{array}{ll}
j & \mbox{if $j\in J$}\, ,\\ j+1 & \mbox{otherwise}\, .
\end{array} \right.
\end{equation}
We are going to construct a locally convex domain $\tilde{R}_{I,J}$ in
 $S$ that ``forces'' the $j$'th strand of a curve homotopic to
 $\gamma_{m,n}$ to cross $E_n$ strictly between 
the points $P_{\ell_j}$ and $P_{r_j}$.  
For each of the bump centers $P_{j}$ there is (on $M^2$) a
geodesic segment $\lambda_{j}$ with endpoints on the
left side of $R_{n}$ which crosses $E_{n}$ to the right
of $P_{j}$, and a geodesic segment $\psi_{j}$  with
endpoints on the right side of $R_{n}$ which crosses $E_{n}$
to the left of $P_{j}$.  Let $\Lambda_{j}$ be the
part of $R_n$ to the left of $\lambda _{j}$, and
$\Psi_{j}$ the part to the right of $\psi_{j}$.  By
construction $\Lambda_{j-1}$ and $\Psi_{j}$ have disjoint
closures and $E_{n}\subset \Lambda_{j}\cup \Psi{_j}$.
Let $\Lambda_{0}$ and $\Psi_{2m+1}$ each be the empty set.
We define the domain $\tilde{R}_{I,J}$ to be $S$ with
the lifts of $\cup_{k\leq \ell_{j}}\Lambda_{k}$ and
$\cup_{k\geq r_{j}}\Psi_{k}$ removed from each
$T_{j}$.  
  
It will be important that these unions are monotonically increasing
(respectively decreasing) with $j$, and that, by definition, 
\begin{equation}
\tilde{R}_{I,J} \cap \tilde{R}_{K,L}=\tilde{R}_
{I\cup K,J\cup L}\, .
\end{equation}

Given $I$ and $J$, subsets of $\{1, \ldots ,2m\}$ again,
we define the face
\begin{equation}
F_{I,J}=\{\epsilon \in [0,1]^{2m} \mid \epsilon_{i} =
0
\mbox{ if $i \in I$ and }\epsilon_{i} =1 \mbox{ if }
i\in J\}\, .
\end{equation}
Note that $F_{I,J}$ is empty unless $I$ and $J$ are
disjoint, and that $F_{I,J}$ is a vertex point  if
they are disjoint with union $\{1, \ldots ,2m\}$.

\begin{Lem} \label{l:l5.3}
There is a continuous map 
$\gamma$ from  $[0,1]^{2m}$ to $H^1(\SS^1, S)$
with the property that if
$\epsilon =(\epsilon_{1}, \ldots ,\epsilon_{2m})\in F_
{I,J}$, then $\gamma
(\epsilon)$ is a
simple closed curve in the interior of $\tilde{R}_{I,J}$, which is
freely homotopic in $S$ to $\tilde{\gamma}_{m,n}$, via
curves whose images in $M^2$ are simple.
\end{Lem}

Note that by the definition of $\tilde{R}_{I,J}$,
if $\epsilon \in F_{I,J}$, the curve $\gamma (\epsilon)
$
described in the lemma will cross the component of
$\tilde{E}_{n}$
containing $Q_{j}$ at a point whose image in $M^2$ lies
(strictly)
between $P_{\ell_{j}}$ and
$P_{r_{j}}$.

\begin{proof}
(of Lemma \ref{l:l5.3}).
The strands are of course ordered not consecutively,
but in the order in which they cross the equator. Note
that $\tilde{R}_{I,J}$ is locally geodesically
convex by
construction. Without spoiling
local convexity we can, by adding an extra geodesic
segment to the bottom (respectively top) of $\lambda$
(respectively $\psi$) assume that the coordinate $z$ is
monotonic along $\lambda$ and $\psi$.  For a fixed
vertex point
$F_{I,J}$ ($I\cup J=\{1,\cdots,2m\})$ of the cube, the
surface
$S$,
(where $\tilde{\gamma}_{m,n}$ lies), can be contracted
along the curves $z=$ constant onto the surface $\tilde
{R}_{I,J}$ keeping $\tilde{R}_{I,J}$ fixed.
The image of $\tilde{\gamma}_{m,n}$ under the
contraction will be $\gamma(\epsilon)$.
Observe that the curve $\gamma(\epsilon)$ agrees with 
$\tilde \gamma_{m,n}$ except where it is pushed aside 
by the curves $\lambda$ and $\psi$.  
Since the sets
$\bigcup_{j\leq k}\Lambda_{j+\epsilon_{k}-1}$
and $\bigcup_{j\geq k}\Psi_{j+\epsilon_{k}}$
are monotonically increasing (respectively decreasing)
with $k$, the image of $\tilde{\gamma}_{m,n}$ under the
contraction, while not necessarily a simple curve, can
be approximated by simple closed curves. We use
linear interpolation along the curves $z=$ constant to
define $\gamma(\epsilon)$ for non-vertex points
$\epsilon$ in the domain $[0,1]^{2m}$.
It is not difficult to see that linear interpolation
will preserve the property that the curves ``move to
the right" with $k$.  Clearly $\tilde \gamma_{m,n}$ 
is not homotopic to a curve in $\partial \tilde R_{I,J}$.  By turning 
on the curve-shortening process momentarily, we can assume that the 
image of $F_{I,J}$ under $\gamma$ consists of simple closed curves 
lying in the interior of $\tilde R_{I,J}$. 
\end{proof}

We next show that the image of the face $F_{I,J}$ under
the map $\gamma$ is nontrivial in the appropriate
   homology group of dimension $2m-|I|-|J|$.  Since
the cube $[0,1]^{2m}$ is contractible, this will have
to be relative homology.  The idea will be to show,
using induction, that $\partial \gamma  F_{I,J}=\gamma \partial  F_{I,J}$ is
nontrivial and to thus conclude that $\gamma F_{I,J}$
is nontrivial.  The nontriviality of $\gamma \partial F_
{I,J}$ will follow from the next lemma.  Let $\Omega_{I,J}$ be the set
of simple closed curves in the interior of $\tilde{R}_{I,J}$ 
which are freely homotopic in $S$ to $\tilde
{\gamma}_{m,n} $ through curves whose images on $M^2$
are simple, and which cross each component of $\tilde{E}
_{n}$ exactly once,
transversely.  Thus:  The $j^{th}$ strand of a curve in
$\Omega_{I,J}$ is constrained to pass to the left of
the point $P_{j}$ if $j\in I$, and to the right of $P_
{j}$ if $j\in J$; it must pass strictly between $P_{j-1}$ and
$P_{j+1}$ in any case.  We can unambiguously
choose a
parameterization for all curves in $\Omega_{I,J}$ by
insisting
that the curves begin on a particular component of
$\tilde{E}_{n}$.
Note that the sets $\Omega_{I,J}$ have the same
intersection
pattern as the faces
$F_{I,J}$ and the domains $R_{I,J}$, namely
\begin{equation}  \label{e:e5.12}
\Omega_{I,J} \cap \Omega_{K,L}=\Omega_{I\cup
K,J\cup L}\, .
\end{equation}
   Also note that a curve in $\Omega_{I,J}$ will cross
the component of
$\tilde{E}_{n}$ containing $Q_{j}$ at a point whose
image in $M^2$ lies
(strictly) between $P_{\ell_{j}}$ and $P_{r_{j}}$.
Finally it will be crucial of course that the sets
$\Omega_{I,J}$ are preserved under the curve shortening
process.  Let
$\Delta_{I,J}$ be the path space where
the boundary of $F_{I,J}$ lies, namely
\begin{equation}
\Delta_{I,J} = \cup_{(I,J)\subset(K,L)}
\Omega_{K,L}\, ,
\end{equation}
where we write $(I,J)\subset(K,L)$ if
$I\subseteq K$ and $J\subseteq L$, but $(I,J)\neq (K,L)$.

\begin{Lem} \label{1:15.4}
The free Abelian group $\langle F_{I,J}\rangle_{\mid
I\mid +\mid J\mid =p}$
   injects into
\begin{equation}
\bigoplus_{|I|+|J|=p}H_{2m-p}
(\Omega_{I,J},\Delta_{I,J})
\end{equation}
and hence injects into
\begin{equation}  \label{e:frst}
H_{2m-p}(\bigcup_{\mid I\mid +\mid J\mid =p}
\Omega_{I,J},\bigcup_{\mid I\mid +\mid J\mid =p}\Delta_
{I,J})\, .
\end{equation}
\end{Lem}

\begin{proof}
Induction (downward) on $p$, using the fact that
(inductively) the boundary $\partial F_{I,J}$ of $F_
{I,J}$ represents a
nontrivial element in $H_{2m-p-1}(\Delta_{I,J})$, and
excision,
together with \eqr{e:e5.12}, for the second injection.
To start the induction we use the fact that the sets
$\Omega_{I,J}$ with $|I|+|J|=2m$ are disjoint by \eqr{e:e5.12}.
\end{proof}

Now let $g_{i}$ be the bumpy metrics close to the given
metric on $M^2$; for simplicity assume also that
distinct closed geodesics have different lengths in the
$g_{i}$ metric.  We are still assuming $m$ and $n$
fixed, with $n$ sufficiently large. If $i$ is
sufficiently large, with the bumpy metric we will (using curve-shortening, 
and the fact that $\eta_n$ has minimal length) have
a locally convex domain $R_{n}$ and geodesics
$\lambda_{j}$ and $\psi_{j}$ with the same pattern
of intersections as in the
limit $i\rightarrow\infty$. All these barriers will
approach the barriers on $M^2$ as $i\rightarrow
\infty$.

\begin{Lem} \label{l:l5.16}
Let $R$ be a locally convex domain, and let $\tau$ be a
geodesic in $R$ with boundary in $\partial R$.  Let
$\gamma$ be a simple closed curve on $R$ intersecting
$\tau$ once, transversely.  Let $ \Omega$ be
the space of (unparameterized) simple closed curves on
the interior of $R$ which are freely homotopic to $\gamma$, and which
intersect $\tau$ once, transversely.
Assume the metric
on $R$ is bumpy, and that different closed geodesics
have different lengths.  Suppose that $a$ and $b$ are
not critical values of the length function, and
that $(a,b)$ contains at most one critical value $c$.
Let $\Omega^{a}$ denote the curves of length $\leq a$.
Then
\begin{equation}  \label{e:e5.18}
H_{k}(\Omega^{b},\Omega^{a})=\left\{ \begin{array}{ll}
                 {\bf Z} & \mbox{if there is a critical
point of index $k$}\, ,\\
0 & \mbox{otherwise}\, .
\end{array}
\right.
\end{equation}
\end{Lem}

\begin{proof}
The curve $\tau$ is used to choose a parameterization
for the curves; we can assume they all start on
$\tau$.  Standard Morse theory arguments, for example
Theorems
4.1 and 4.2 on p. 34--35 in Chang's book \cite{Ch}, which use the
gradient of the energy function on the space of $H^1$
curves apply in this context as well, given the
previously mentioned properties of the curve-shortening
flow, i.e., \eqr{e:kappainfty}.  In fact the
computation of the local critical group needs no
alteration since a sufficiently small neighborhood of a
simple closed geodesic will consist entirely of simple
curves.
\end{proof}

Again let $m$ and $n$ be fixed with $n$ sufficiently
large.  In the $g_{i}$ metric for $i$ sufficiently large
there will be locally convex domains $R_{I,J}$ and path
spaces $\Omega_{I,J}$.  With these hypotheses we have:

\begin{Lem} \label{l:l5.4}
Each $\Omega_{I,J}$ contains a (``minimax") geodesic of
index $2m-p$, where $|I|+|J|=p$.
\end{Lem}

\begin{proof}
Let $\delta$ be small but positive.  We use downward
induction on $p$.  The
idea is to push
down each cell $\gamma F_{I,J}$ using the curve-shortening flow, leaving
its (already pushed down) boundary fixed. It seems
convenient however to work
with relative homology groups.  The cells $\gamma F_{I,J}$ can be converted 
into singular chains using the shuffle homeomorphism (\cite{GgHa}, p. 268).  
All singular chains will have $\ZZ_2$ coefficients.  We have the boundary 
relations
\begin{equation}  \label{e:st}
\partial \gamma F_{I,J}
=\sum_{(I,J)\subset (K,L)\, ,\,|I|+|J|+1=|K|+|L|}\gamma F_{K,L}\, .
\end{equation}
We will inductively
alter (``pushdown'') the cells $\gamma
F_{I,J}$ by replacing each by a homologous (in $\Omega_
{I,J}$) chain in
$\Omega_{I,J}$.  If we add $\partial\,\tau$ to $\gamma F_{K,L}$, 
for a chain $\tau$ in $\Omega_{K,L}$, then we must add $\tau$ 
to each $\gamma F_{I,J}$ with $(I,J)\subset (K,L)$ and 
$|I|+|J|+1=|K|+|L|$; this way the boundary relations \eqr{e:st} 
will be maintained.

The lemma is clear from the Barrier Principle when
$p=2m$ since
$\Omega_{I,J}$ is
nonempty.  For each $I,J$ with $p=2m$, the minimum
value $a_{I,J}$ with $H_{0}^{a_{I,J}}(\Omega_{I,J})
\not=0$ is
clearly a critical value of the length function
corresponding to a critical point of index $0$. We can
(using
curve-shortening) assume that $\gamma F_{I,J}$ lies in
$\Omega_{I,J}^{a_{I,J}+\delta}$.

   Now fix $I,J$ with $|I|+|J|=p$.
Assume that for each $(K,L)\supset (I,J)$ we have
$\gamma F_{K,L}\subset \Omega_{K,L}^{a_{K,L}+\delta}$
and
assume that $a_{K,L}-\delta > a_{M,N}+\delta$ if $(K,L)
\subset (M,N)$.
Assume that in each $\Omega_{K,L}\setminus \Delta_{K,L}$ with 
$(I,J)\subset (K,L)$, we have a critical point of index $2\,m-|K|-|L|$ 
whose critical value $a_{K,L}$ is the infimum $a$ for which the class 
of $\gamma \partial F_{K,L}$ vanishes in $H_{\ast}(\Omega_{K,L}^a)$.  
This means that the image of $\gamma \partial F_{K,L}$ is nontrivial 
in $H_{\ast}(\Omega_{K,L}^{a_{K,L}-\delta})$, and that 
$H_{\ast}(\Omega_{K,L}^{a_{K,L}+\delta},\Omega_{K,L}^{a_{K,L}-\delta})
\ne 0$.  The assumption that the critical point with critical value 
$a_{K,L}$ lies in $\Omega_{K,L}\setminus \Delta_{K,L}$, together with 
the assumption that the different critical points have different 
critical values, and \eqr{e:e5.12}, implies that $a_{M,N}\ne a_{K,L}$ 
if $(M,N)\ne (K,L)$.
The
lemma will follow from
Lemma \ref{l:l5.16} once we establish that
\begin{equation}  \label{e:imagenontr}
\mbox{the image of }\gamma \partial F_{I,J}
\mbox{ is nontrivial in }H_{\ast}(\Omega_{I,J}^
{a+\delta})\, ,
\end{equation}
where $a=\max_{(I,J)\subset (K,L)}a_{K,L}$, and 
\begin{equation}  \label{e:stst}
\mbox{the image of }\gamma \partial F_{I,J}
\mbox{ is nontrivial in }H_{\ast}(\Delta_{I,J})\, .
\end{equation}

In order to establish \eqr{e:imagenontr} we argue as
follows.
Suppose
$a=\mbox{(say)}a_{M,N}$.  Line up the long exact
sequences for the pairs
$(\Delta_{I,J}^{a+\delta},\Delta_{I,J}^{a-\delta})$ and
$(\Omega_{I,J}^{a+\delta},\Omega_{I,J}^{a-\delta})$.
The class of $\gamma F_{M,N}$ in 
$H_{\ast}(\Delta_{I,J}^{a+\delta},\Delta_{I,J}^{a-
\delta})$ is the same as
that of $\gamma \partial F_{I,J}$, which comes from
$H_{\ast}(\Delta_{I,J}^{a+\delta})$.  The further image
in
$H_{\ast}(\Omega_{I,J}^{a+\delta},\Omega_{I,J}^{a-
\delta})$ is nonzero by
Lemma \ref{l:l5.16}, and the inductive hypothesis.  
Thus the image of $\gamma
\partial F_
{I,J}$ in
$H_{\ast}(\Omega_{I,J}^{a+\delta})$ is nonzero.  
We get \eqr{e:stst} from \eqr{e:frst} by looking at the image of 
\begin{equation}
H_{\ast}(\Delta_{I,J})=H_{\ast}(\bigcup_{(I,J)\subset (K,L)}\Omega_{K,L})
\end{equation}
in 
\begin{equation}
H_{\ast}(\bigcup_{|K|+|L|=p+1}\Omega_{K,L},
\bigcup_{|K|+|L|=p+1}\Delta_{K,L})\, .
\end{equation}
On the
other hand since $\partial \gamma F_{I,J}=\gamma
\partial F_{I,J}$, the image of $\gamma \partial F_{I,J}
$ in $H_{\ast}(\Omega_{I,J})$ is $0$; hence for some
$b>a$, $H_{2m-p}(\Omega_{I,J}^{b},\Omega_{I,J}^{a})\ne 0$,
so that we can use \eqr{e:e5.18} to 
get a critical point in $\Omega_{I,J}$ of index $2\,m-p$.  
That the critical point does not lie in $\Delta_{I,J}$ follows 
from \eqr{e:stst}.
\end{proof}

\begin{proof}
(of Theorem \ref{t:counterex}).
For any $m\geq 1$, if $n$ is sufficiently large and
$p\leq 2m$, , we get a sequence $\{\sigma_{i}\}$ of
simple closed curves with $\sigma_{i}$ a geodesic in the
$g_{i}$ metric, with index $2m-p$.  A subsequence will
converge to a closed simple geodesic $\sigma$ with
index $2m-p-2$, $2m-p-1$, or $2m-p$ and which is freely
homotopic in
$S$ to $\tilde{\gamma}_{m,n}$.  In fact
since for different $n$ these are different free
homotopy classes in the cylinder-minus-$3$-holes, for
each $k=2m-p$ there are an infinite number of simple
closed geodesics of index $k-2, k-1,$ or $k$.
\end{proof}

\begin{Rem}   \label{r:genright}
The metric $M_0^2$ described above is not bumpy, as the 
``middle sphere section'' contains a neighborhood of a 
great circle on the standard sphere.  However there is a 
lot of inessential symmetry in the construction.  The curves 
$\eta_n$ and $\gamma_{m,n}$ will exists for all metrics in a 
neighborhood of $M_0^2$.  In order to get simple closed geodesics of 
arbitrary high Morse index, all we really need is what comes from 
Lemma \ref{l:l5.2}:  For each $m$, if $n$ is sufficiently large we need 
$2\,m$ points $P_1,\cdots,P_{2m}$ on $E_n$ and, for each $j<2m$, a 
geodesic $\lambda$ on $R_n$ with endpoints on the left side of $R_n$, 
and crossing $E_n$ once between $P_j$ and $P_{j+1}$; and a geodesic 
$\psi$ on $R_n$ with endpoints  on the right side of $R_n$, and 
crossing $E_n$ once between $P_j$ and $P_{j+1}$.  While we do not know 
how to prove that there is a bumpy metric with this property, we see no 
reason why  one should not exist.  The construction requires that some 
$\underline{\text{\it noncompact}}$ simple geodesic is a limit of unstable 
geodesics; for this bumpyness of the metric seems irrelevant.
\end{Rem}

\appendix

\section{Convexity of a neighborhood of a strictly stable geodesic}
\label{s:convex}

Let $\gamma\subset M$ be a simple closed geodesic and let $(s,\theta)$ be
Fermi coordinates in a neighborhood of $\gamma=\gamma (\theta)$.
In these coordinates
the metric can be
written as $ds^2+f^2(s,\theta)\,d\theta^2$.  For $\alpha>1$ and
$\phi=\phi (\theta)>0$ set
$F(s,\theta)=s^{\alpha}\,\phi^{-\alpha}(\theta)$, then
\begin{equation}
\nabla F
=\alpha s^{\alpha-1}\phi^{-\alpha}\frac{\partial}{\partial s}
-\alpha s^{\alpha}\phi_{\theta}\phi^{-\alpha-1}f^{-2}
\frac{\partial}{\partial \theta}\, ,
\end{equation}
\begin{equation}
\langle \nabla_{\frac{\partial}{\partial s}}\nabla F,
\frac{\partial}{\partial s}\rangle
=\alpha (\alpha-1)s^{\alpha-2}\phi^{-\alpha}\, ,
\end{equation}
\begin{equation}
\langle \nabla_{\frac{\partial}{\partial s}}\nabla F,
\frac{\partial}{\partial \theta}\rangle
=\alpha\,\frac{\phi_{\theta}}{\phi^{\alpha+1}}\, s^{\alpha-1}
\,\left(s \frac{f_s}{f}
-\alpha\right)\, ,
\end{equation}
\begin{equation}
\langle \nabla_{\frac{\partial}{\partial \theta}}\nabla F,
\frac{\partial}{\partial \theta}\rangle
=\alpha s^{\alpha-1}\frac{f\,f_s}{\phi^{\alpha}}
-\alpha s^{\alpha}\frac{\phi_{\theta,\theta}}{\phi^{\alpha+1}}
+\alpha (\alpha+1) s^{\alpha}\frac{\phi^2_{\theta}}{\phi^{\alpha+2}}
+\alpha 
s^{\alpha}\frac{\phi_{\theta}}{\phi^{\alpha+1}}\frac{f_{\theta}}{f}\, .
\end{equation}
Since $K\,f=-f_{s,s}$ we get by Taylor expansion that
\begin{equation}
\langle \nabla_{\frac{\partial}{\partial \theta}}\nabla F,
\frac{\partial}{\partial \theta}\rangle
=-\alpha s^{\alpha}\,\frac{L_{\gamma}\phi}{\phi^{\alpha+1}}
+\alpha (\alpha+1) s^{\alpha}\frac{\phi^2_{\theta}}{\phi^{\alpha+2}}
+\frac{\alpha}{\phi^{\alpha}} o(s^{\alpha})\, .
\end{equation}
   From this it follows easily that if $\gamma$ is strictly stable and
$\phi$ is a positive eigenfunction of $L_{\gamma}$ with eigenvalue
$\lambda_1(L_{\gamma})>0$,
then for $\alpha>1$ sufficiently large
$F$ is convex in a neighborhood $T$ of $\gamma$
and strictly convex in $T\setminus \{\gamma\}$.

\vskip2mm
Recall that we equip the space of $C^{\infty}$ metrics on a closed surface
$M^2$ with the
$C^{\infty}$-topology and we write $g_i\to g$ if $|g-g_i|_{C^{\infty}}\to
0$.

An easy consequence of the existence of $F$ is:

\begin{Lem}    \label{l:exF}
Let $M^2$ be a closed surface with a metric $g$ and suppose that 
$\gamma\subset M$ is a simple closed strictly stable geodesic.   
Then there exists $\epsilon>0$ and an open neighborhood $T$ of 
$\gamma$ such that if $\tilde g$ is a metric on $M$ with 
$|g-\tilde g|_{C^{\infty}}<\epsilon$, then there exists a 
unique closed simple geodesic (in the metric $\tilde g$) 
$\tilde \gamma\subset T$.  Moreover, $\epsilon$ and $T$ 
can be chosen so that $\tilde \gamma$ is strictly stable. 
\end{Lem}

\begin{proof}
First using the existence of $F_{\gamma}$ it follows easily that given a 
metric $\tilde g$ sufficiently close to $g$ there must be a simple 
closed strictly stable geodesic $\tilde \gamma$.  Now using 
$F_{\tilde \gamma}$ it follows easily that $\tilde \gamma$ is unique.  
\end{proof}

\section{Convergence of metrics and geodesics}

\begin{Lem}  \label{l:vary}
Let $M^2$ be as above with a bumpy metric $g$.  For each $L>0$,
there exists at most finitely many closed geodesics of length $<L$.
Moreover, if $L$ is not equal to the
length of any closed geodesic in $g$, then in a
neighborhood of $g$ each metric has precisely as many (simple)
closed stable geodesics of length $<L$ as $g$.  Finally, if
$g_i\to g$ and $\{\gamma_{i,k}\}$, $\{\gamma_k\}$
are the (simple)
closed stable geodesics in
$g_i$, $g$, respectively, of length $<L$,
then $\gamma_{i,k}\to \gamma_k$ for $i\to \infty$ and
each $k$.
\end{Lem}

\begin{proof}
If there were an infinite sequence of
such (simple) closed geodesics, then it would follow that a
subsequence would converge to a closed
geodesic with a nontrivial Jacobi field contradicting that the
metric is bumpy.  In fact it is easy to see
(by locally going to a finite cover) that the assumption that
the geodesics are simple is not needed for this
or anything else in this lemma.

It follows easily from the existence of the convex function $F$ from
Appendix \ref{s:convex} that for each metric in a neighborhood of $g$ there
are at least as many (simple)
closed stable geodesics of lengths at most $L$ as in $g$.
That there are not more (and the last claim) follows from a simple
convergence argument together with the assumption that the metric is bumpy.
\end{proof}

\end{document}